\documentclass[11pt,a4paper]{article}
\usepackage{fullpage}
\usepackage{graphicx}
\usepackage{amsmath}
\usepackage{amsfonts}
\usepackage{amssymb}
\newtheorem{theorem}{Theorem}

\newtheorem{claim}[theorem]{Claim}

\newtheorem{corollary}[theorem]{Corollary}

\newtheorem{definition}[theorem]{Definition}

\newtheorem{lemma}[theorem]{Lemma}

\newtheorem{proposition}[theorem]{Proposition}

\newenvironment{proof}[1][Proof]{\textbf{#1.} }{\ \rule{0.5em}{0.5em}}

\newcommand{\R}{\ensuremath{\mathbb{R}}}
\newcommand{\Z}{\ensuremath{\mathbb{Z}}}
\newcommand{\T}{\ensuremath{\mathbb{T}}}

\newcommand{\A}{\mathcal{A}}
\newcommand{\Q}{\ensuremath{\mathbb{Q}}}
\newcommand{\eps}{\varepsilon}
\renewcommand{\P}{\mathbb{P}}
\newcommand{\E}{\mathbb{E}}

\newcommand{\nchoosek}[2]{\begin{pmatrix}#1\\#2\end{pmatrix}}
\newcommand{\ra}{\ensuremath{\rightarrow}}

\newcommand{\ol}{\overline}

\renewcommand{\S}{\mathcal{S}}
\DeclareMathOperator{\Maj}{Maj}
\DeclareMathOperator{\var}{Var}
\DeclareMathOperator{\bin}{Bin}
\DeclareMathOperator{\Tribes}{Tribes}

\begin{document}
\thispagestyle{empty}

\title{On K-wise Independent Distributions and Boolean Functions}

\author{{Itai Benjamini \thanks{Weizmann Institute, Rehovot, 76100,
Israel. itai.benjamini@weizmann.ac.il}} \quad {Ori Gurel-Gurevich
\thanks{Weizmann Institute, Rehovot, 76100, Israel.
ori.gurel-gurevich@weizmann.ac.il}} \quad {Ron Peled \thanks{UC Berkeley. peledron@stat.berkeley.edu}}}

\maketitle

\begin{abstract}
We pursue a systematic study of the following problem. Let
$f:\{0,1\}^n \rightarrow \{0,1\}$ be a (usually monotone) boolean
function whose behaviour is well understood when the input bits are
identically independently distributed. What can be said about the
behaviour of the function when the input bits are not completely
independent, but only $k$-wise independent, i.e. every subset of $k$
bits is independent? more precisely, how high should k be so that
any k-wise independent distribution "fools" the function, i.e.
causes it to behave nearly the same as when the bits are completely
independent?

In this paper, we are mainly interested in asymptotic results about
monotone functions which exhibit sharp thresholds, i.e. there is a
critical probability, $p_c$, such that $P(f=1)$ under the completely
independent distribution with marginal $p$, makes a sharp
transition, from being close to 0 to being close to 1, in the
vicinity of $p_c$. For such (sequences of) functions we define 2
notions of "fooling": $K_1$ is the independence needed in order to
force the existence of the sharp threshold (which must then be at
$p_c$). $K_2$ is the independence needed to "fool" the function at
$p_c$.

In order to answer these questions, we explore the extremal
properties of $k$-wise independent distributions and provide ways of
constructing such distributions. These constructions are connected
to linear error correcting codes.

We also utilize duality theory and show that for the function $f$ to
behave (almost) the same under all $k$-wise independent inputs is
equivalent to the function $f$ being well approximated by a real
polynomial in a certain fashion. This type of approximation is
stronger than approximation in $L_1$.

We analyze several well known boolean functions (including AND,
Majority, Tribes and Percolation among others), some of which turn
out to have surprising properties with respect to these questions.

In some of our results we use tools from the theory of the classical
moment problem, seemingly for the first time in this subject, to
shed light on these questions.

\end{abstract}

\thispagestyle{empty}

\newpage
\setcounter{page}{1}

\section{Introduction}
Let $f:\{0,1\}^n \rightarrow \{0,1\}$ be a boolean function whose
behaviour is well understood when the input bits are independent and
identically distributed, with probability $p$ for each bit to be 1.
As an example we may consider the majority function, $\Maj$, whose
output is the bit which occurs more in the input (suppose that $n$
is odd). When $p=1/2$ we know that the output is also distributed
uniformly. When $p<1/2$ the output tends to be 0. More precisely, if
$p<1/2$ is constant, the probability of $\Maj=1$ decays
exponentially fast with $n$.

Suppose, however, that the input bits are not truly IID . For
example, they might be the result of a derandomization procedure. A
reasonable, but weaker assumption would be that the probability of
each bit to be 1 is still $p$, and that they are $k$-wise
independent, i.e. the distribution of any $k$ of the bits is
independent.

Under this assumption, what can be said about the distribution of
$f$? For fixed $p$, which $k$ (as a function of $n$) is enough to
guarantee the same asymptotic behaviour? Majority turns out to be
relatively easy to analyze: $k=2$ is enough to guarantee that for
any fixed $p<1/2$, the probability of $\Maj=1$ tends to 0 (though
only polynomially fast), while for $p=1/2$, we have $P(\Maj=1)$
guaranteed to tend to $1/2$ if and only if $k=\omega(1)$
("Guarantee" here means that $\Maj$ behaves as prescribed under any
$k$-wise independent distribution). In fact, for $p=1/2$ we have
more precise results, that $|P(\Maj=1) - 1/2| \le O(1/\sqrt{k})$
under any $k$-wise independent distribution. As can been seen, the
$k$ needed to "fool" majority at $p\ne p_c$ (which we denote $K_1$)
is much smaller then the $k$ needed to "fool" majority at $p_c$ (
which we denote $K_2$). This phenomenon is shared by the other
functions we explore, and we provide a partial explanation. Other
functions exhibit much more complex behaviour and the required
analysis is accordingly complex. We pursue a systematic study of the
above question.

$k$-wise independent distributions are often used in computer
science for derandomization of algorithms. This was initiated by the
papers \cite{ABI86}, \cite{CGHFRS85}, \cite{KW85}, \cite{L85} and
further developed in \cite{BR91}, \cite{MNN94}, \cite{L93},
\cite{S92}, \cite{KM93}, \cite{KM94} and others (see \cite{LW95} for
a survey). For derandomization one checks that the algorithm still
behaves (about) the same on a particular $k$-wise independent input
as in the completely independent case. The question we ask is of the
same flavor, for a given boolean function $f$, we ask how much
independence is required for it to behave about the same on all
$k$-wise independent inputs (including the completely independent
one).

%
%In order to derandomize a certain algorithm using a certain
%$k$-wise independent distribution one needs to check that the
%algorithm still behaves similarly on this $k$-wise independent
%input as it did for the fully random input. The question we ask is
%of the same flavor, for a given boolean function $f$, the
%algorithm, we ask how much independence is required for it to
%behave about the same on all $k$-wise independent inputs
%(including the fully random one).

Typically, $k$-wise independent distributions are constructed by
sampling a uniform point of a small sample space, which is usually
also a linear subspace (\cite{J74}, \cite{HSS99}, \cite{MS77}). In
this work, like in the works of \cite{KM93}, \cite{KM94}, we do not
impose this restriction and consider general $k$-wise independent
distributions. Still, our work is of interest even for the reader
only interested in the more restrictive model since, on the one
hand, anything we show is impossible would still be impossible in
that model and on the other hand almost all of our constructions are
of the linear subspace type. Interestingly, in section
\ref{max_prob_all_bits_1_sec} we give an example where the general
and more restrictive case give asymptotically different results,
i.e., that the general distribution case is richer, not just up to
constants, in what can be achieved with it.

The tools we use include the duality of linear programming, in
section \ref{duality_sec}, used to show an \emph{equivalence}
between our question and the question of approximating the function
$f$ by a real polynomial in a certain "sandwich $L_1$" approximation
(stronger than ordinary $L_1$ approximation). This connects our
results to the subject of approximation of boolean functions, used
for example in learning theory (e.g. \cite{LMN89}, \cite{NS92},
\cite{ABFR91}, \cite{OS03}).

In section \ref{dist_determined_by_moments_sec} we recall a theorem
about weak convergence of distributions, later used to give sharp
bounds on $K_2$ very easily. In section \ref{TCMP} we introduce a
tool from the Theory of the Classical Moment Problem (TCMP),
seemingly for the first time in this context. In sections
\ref{max_prob_all_bits_1_sec} and \ref{min_prob_all_bits_1_sec} we
use it to prove bounds on the maximal and minimal probabilities of
all bits to be $1$ under a $k$-wise independent distribution, in a
simple way. We then observe that if $p=\frac{1}{q}$ for a
prime-power $q$, then an upper bound on this maximal probability
translates to a lower bound on the size of a symmetric sample space
for $k$-wise independent $GF(q)$-valued random variables, we apply
our upper bound to obtain new lower bounds for such sample spaces.
For the binary case $q=2$ our bound equals the well-known bound of
\cite{ABI86}, \cite{CGHFRS85}.

%By the equivalence if $f$ behaves about the same on all $k$-wise
%independent inputs then it can be well approximated by degree $k$
%polynomials (in our sense) and if it changes its behaviour on some
%$k$-wise independent input then it cannot.

%In sections \ref{max_prob_all_bits_1_sec} and
%\ref{min_prob_all_bits_1_sec} we investigate the maximal and minimal
%probabilities of all bits to be $1$ under a $k$-wise independent
%distribution. We present new bounds for these extremal
%probabilities. We then observe that if $p=\frac{1}{q}$ for a
%prime-power $q$, then an upper bound on this maximal probability
%translates to a lower bound on the size of a symmetric sample space
%for $k$-wise independent $GF(q)$-valued random variables, we apply
%our upper bound to obtain new lower bounds for such sample spaces.
%For the binary case $q=2$ our bound equals the well-known bound of
%\cite{ABI86}, \cite{CGHFRS85}. In this section we use a tool
%apparently not used before in this context, that of the Theory of
%the Classical Moment Problem ("TCMP"). This tool allows us to obtain
%some of our results very simply.

In section \ref{boolean_functions_sec} we explore $K_1$ and $K_2$
for various boolean functions, and also prove some general theorems.
In section \ref{percolation_section} we present a novel construction
of a distribution (of the linear subspace type) designed to change
the behaviour of a particular function. We use a variation of the
$(u\ |\ u+v)$ construction of error-correcting codes \cite{MS77} and
we would like to emphasize the technique used there. We think there
is a shortage of ways to construct $k$-wise independent
distributions with specified properties and that this technique will
be useful for changing the behaviour of other functions as well.

The approach in this paper is a little different than that usually
taken in pseudo-random generators (see \cite{NW94}). There one seeks
a distribution under which all functions from a certain complexity
class behave the same as on fully independent bits. In contrast, we
start with a function $f$ and wish to show that it behaves the same
on all $k$-wise independent inputs. Still, one may expect this to
hold if the function $f$ is "simple enough". Indeed, a conjecture of
Linial and Nisan \cite{LN} makes this precise when $f$ is a function
from the class AC0. In section \ref{AC0_functions} we recall the
precise conjecture and make some modest progress towards confirming
it.

There are other notions of "simple functions". Another such notion
is that the function be noise stable \cite{BKS}, i.e., having most
of its Fourier mass on constant level coefficients. In section
\ref{fourier_sec} we show a connection between the Fourier spectrum
and the behaviour on $k$-wise independent inputs, but surprisingly
show that a noise stable function can behave very differently on
$k$-wise independent inputs than on fully independent inputs even
when $k$ grows fast with $n$.

There is also a lot of interest in \emph{almost} k-wise independent
distributions
(\cite{NN93},\cite{AGHP92},\cite{ABNNR92},\cite{AMN90},\cite{CRS94},\cite{EGLNV92}),
though our questions can equally be formulated for that case, in our
work we concentrate only on perfect $k$-wise independence, this is
both because it seems the analysis is simpler for perfect $k$-wise
independence and they could serve as a starting point for further
research and because we think the perfect $k$-wise independent case
is interesting on its own.

%Table of contents:
%\section{Basic definitions and properties}
%\section{General tools}
%\section{Boolean functions}
%\subsection{Majority}
%\subsection{AC0 functions}
%\subsubsection{Tribes}
%\subsubsection{General theorem}
%\subsubsection{A more complex example}
%\subsection{Composition of functions}
%\subsection{Fourier transform and K_2}
%\section{Probabilities that all bits are 1}

\section{Basic definitions and properties} \label{basic_sec}
We begin with a definition
\begin{definition}
Let $\A(n,k,p)$ be the set of all $k$-wise independent
distributions $\Q$ on $n$ bits $(X_1,\ldots, X_n)$ with
$\Q(X_i=1)=p$ for all $i$.
\end{definition}
Also denote by $\P_p$ the fully independent distribution on $n$
bits, each with probability $p$ to be $1$.

In most of the sequel we will be concerned with understanding
\begin{equation} \label{extremal_values_for_functions}
\max_{\Q\in\A(n,k,p)} \Q(f=1) \qquad \text{ and } \qquad
\min_{\Q\in\A(n,k,p)} \Q(f=1)
\end{equation}
for a boolean function $f:\{0,1\}^n \rightarrow \{0,1\}$ and given
$n,k$ and $p$. We first note that $\A(n,k,p)$ is a convex set,
since, if the distribution on a subset of $k$ bits is independent
with marginal $p$ in both $\Q_1$ and $\Q_2$ then it is so also in
$\alpha \Q_1 + (1-\alpha) \Q_2$.

%the trivial
%\begin{lemma}
%For every $\Q_1,\Q_2\in\A(n,k,p)$ and $0\le\alpha\le1$ we have
%$\Q:=\alpha \Q_1 + (1-\alpha) \Q_2 \in \A(n,k,p)$.
%\end{lemma}
%\begin{proof}
%Fix a subset of $m\le k$ of the bits. Let $A$ be the event that
%the bits in this subset are all $1$, then $\Q_1(A)=\Q_2(A)=p^m$.
%Hence also $\Q(A)=p^m$. This implies that $\Q\in\A(n,k,p)$.
%\end{proof}
%
%That is, the set $\A(n,k,p)$ is convex.

This implies that the extremal values in
\eqref{extremal_values_for_functions} are attained at extreme points
of $\A(n,k,p)$, hence if we could only find all these extreme points
we could then find the values \eqref{extremal_values_for_functions}
for all $f$. Unfortunately saying anything about these extreme
points appears to be very difficult and so in the sequel we will
need to resort to special methods for each function $f$ considered.

For later reference, we identify the two extreme points of
$\A(n,n-1,\frac{1}{2})$. XOR0 is the distribution on $(X_1,\ldots,
X_n)$ having $\{X_i\}_{i=1}^{n-1}$ IID and
$X_n\equiv\sum_{i=1}^{n-1} X_i$ mod 2, and XOR1 is the same with
$X_n\equiv1+\sum_{i=1}^{n-1} X_i$ mod 2.

We next define precisely what we mean by "$k$ large enough so that
$f$ behaves on all $k$-wise independent inputs the same as on the
fully independent input".
\begin{definition}
$\eps^f(k,p)=\max_{\Q\in\A(n,k,p)} \Q(f=1) - \min_{\Q\in\A(n,k,p)}
\Q(f=1)$

$k^f(\eps,p)$ is the minimal $k$ such that $\eps^f(k,p)<\eps$.
\end{definition}

We will be mostly interested in asymptotic (in $n$) results. Let
$f_n:\{0,1\}^n \rightarrow \{0,1\}$ be a sequence of monotone
boolean functions. Assume that the functions have a sharp
threshold, i.e. there is a $p_c$ such that $\lim_{n\ra \infty}
\P_p(f=1)$ is 0 if $p<p_c$, 1 if $p>p_c$.

%We usually also require the limit to exist when $p=p_c$, and denote
%it by $\alpha$ (*** Do we ever use this line? ***) .

For example, any sequence of balanced monotone transitive
functions has a sharp threshold, as is proved by Friedgut and
Kalai \cite{FK}.

\begin{definition}
$K_1$ is the class of functions $k(n)$, such that $\eps(k,p)\ra 0$
for any $p\ne p_c$.

$K_2$ is the class of functions $k(n)$, such that $\eps(k,p_c)\ra
0$.
\end{definition}

In other words, $K_1$-wise independence is enough to guarantee the
existence of sharp threshold (which is then necessarily at $p_c$),
while $K_2$-wise independence is enough to guarantee that $f$
behaves as if the bits were completely independent, when $p=p_c$.

Notice that while $K_1$ and $K_2$ are classes of functions, we
occasionally abuse the formal notation, and write, as above,
$K_1$-wise independence. Similarly, we write $K_1>k(n)$ to
indicate that $k(n)$ does not belong to $K_1$, or $K_2<\omega(1)$
to indicate $K_2\supset\omega(1)$, etc.

It is not a-priori clear whether $K_1 \le K_2$ or vice versa (or
neither). Consult the appendix for a partial result. In all the
examples we encountered $K_2$ is at least $\omega(K_1)$.

\section{General tools}

In this section we discuss some general tools for finding $K_1$ and
$K_2$ as defined in the previous section.

\subsection{Duality - Approximation by polynomials}
\label{duality_sec}

We note that the values \eqref{extremal_values_for_functions} are
the solution to a simple linear program. What is the dual of this
program? We observe
\begin{proposition}
For any $f:\{0,1\}^n \rightarrow \{0,1\}$, any $k$ and any
$0<p<1$.
\begin{equation}
\begin{split}
\max_{\Q\in\A(n,k,p)} \Q(f=1) &= \min_{P\in P_k^+(f)} \E_{\P_p}
P(X_1,\ldots, X_n)\\
\min_{\Q\in\A(n,k,p)} \Q(f=1) &= \max_{P\in P_k^-(f)} \E_{\P_p}
P(X_1,\ldots, X_n)
\end{split}
\end{equation}
where $P_k^+(f)$ is the set of all real polynomials $P:\R^n\to\R$
of degree not more than $k$ satisfying $P\ge f$ on all points of
the boolean cube. $P_k^-$ is defined analogously with $P\le f$.
\end{proposition}
The proof is simple using linear programming duality. We deduce that
$\eps^f(k,p)< \eps$ is equivalent to having two polynomials $P^+\ge
f$ and $P^-\le f$ of degree not more than $k$ with $\E_{\P_p}(P^+ -
P^-)< \eps$. We call this type of approximation of $f$ a
\emph{"sandwich $L_1$" approximation}. In section \ref{fourier_sec}
we show that it is strictly stronger than $L_1$ approximation (by
real polynomials of degree not more than $k$). Whether it is
stronger than $L_2$ approximation is one of our main open questions.

\subsection{Distributions determined by their moments}
\label{dist_determined_by_moments_sec}
\begin{definition}
We say that a real random variable $X$ has \emph{distribution
determined by its moments} if any random variable $Y$ satisfying
$\E X^m=\E Y^m$ for all integer $m\ge 1$ has the same distribution
as $X$.
\end{definition}
We shall often use the following principle
\begin{proposition} \label{conv_of_moments_prop}
Suppose a sequence of RV's $\{X_n\}_n$ satisfies for all $m$, $\E
X_n^m \to \E X^m$ for some RV $X$ whose distribution is determined
by its moments. Then $X_n\to X$ in the weak sense.
\end{proposition}
For the proof, see \cite{D96}, section 2.3 . We remark that a
distribution is determined by its moments whenever these do not grow
too fast. The best criterion is called Carleman's condition (see
\cite{D96}). But for our purposes it will mostly be enough to know
that the Normal and Poisson distributions are determined by their
moments.

\subsection{Bounds from the classical moment problem} \label{TCMP}
Given a real sequence $\S := \{s_m\}_{m=0}^k$, with $k$ even and
$s_0=1$ , let
%\footnote{this last condition is convenient for us in
%order to use probabilistic notation, but it is not necessary for the
%results of the classical moment problem.}
\begin{equation}
\A_{\S} = \{\Q\ |\ \text{$\Q$ a probability distribution on $\R$,
$s_m = \E_{Q} (X^m)$ for $0\le m\le k$}\}
\end{equation}
be all probability distributions with these first $k$ moments ($X$
is a random variable distributed according to $\Q$). In the theory
of the classical moment problem \cite{Ak65}, \cite{KN77}, based on
$\S$ a certain sequence of real functions $\rho_m$ is defined and
the following theorem is proved
\begin{theorem} \label{maximal_mass_difference_thm} \cite[2.5.2 and 2.5.4]{Ak65}
For any $x$ and any $\Q_1,\Q_2\in\A_{\S}$
\begin{equation}
|\Q_1(X\le x) - \Q_2(X<x)|\le \rho_{\frac{k}{2}}(x)
\end{equation}
and in particular by taking $\Q_1=\Q_2$ we get $\max_{\Q\in\A_{\S}} \Q(X=x)\le \rho_{\frac{k}{2}}(x)$.
\end{theorem}
For brevity we do not give the general definitions of $\rho_m$
here but differ them to the appendix. In the case of interest for us $s_m:=\E(X^m)$ when
$X\sim\bin(n,p)$ and then $\rho_m(x):=(\sum_{j=0}^{m}
P_j^2(x))^{-1}$ where $\{P_j\}_j$ are the (normalized) Krawtchouk
polynomials (see \cite{Sz75}). These polynomials are very well
known and from them we easily deduce the following (see appendix
for a proof)
\begin{align}
\rho_m(n) &= \frac{p^n}{\P(\bin(n,1-p)\le m)}&&\label{rho_last_value}\\
\rho_m(\frac{n}{2}) &\le \frac{2}{\sqrt{m}} &&\text{for
$p=\frac{1}{2}$, even $n$ and even $m\le \frac{n}{2}$}\label{rho_middle_value_bound}
\end{align}
In many cases, the theory also has constructions achieving the bound
of theorem \ref{maximal_mass_difference_thm}. However, these are not
necessarily supported by the integers, which we require. It might be
that they can be suitably modified to give sharp results in our
cases.

%It is possible that a modification of these constructions can remedy
%this.
%
%, but we could not use these since in the cases we needed we
%required the support of the distribution to be on integer points. It
%is possible, however, that a modification of these constructions can
%yield a distribution on integer points, this would be very useful to
%show the sharpness of the bounds in the cases we use.

\section{Boolean functions} \label{boolean_functions_sec}
In this section we investigate $K_1$ and $K_2$ for several boolean
functions, and also present some general theorems. We start with a
simple, but already non-trivial example.

\subsection{Majority}
Let $Maj_n$ be the majority function on $n$ bits (for odd $n$).
Let $S_n=\sum_{i=1}^n x_i$, where $x_i$ are the input bits. Let
$\overline{S_n}=(2 S_n - n)/ \sqrt{n}$. The central limit theorem
implies that under $\P_{1/2}$, $\ol{S_n}\ra N(0,1)$. Identifying
$K_1$ is easy

\begin{theorem} \label{maj_k_1}
$K_1(\Maj)=2$.
\end{theorem}

%\begin{proof}
\noindent\textbf{Proof.} Obviously, $k=1$ is not in $K_1$. However,
for $\Q\in\A(n,2,p)$ we have $\E_\Q(S_n)=np$ and
$\var_\Q(S_n)=np(1-p)$. If, WLOG, $p<1/2$ then by Chebyshev's
inequality

$$\Q(S_n > n/2) \le \Q((S_n -np) > n(1/2-p)) \le
\frac{np(1-p)}{(n(1/2-p))^2} = O(\frac{1}{n}) \ra 0 \quad
\rule{0.5em}{0.5em}$$
% \ \ \blacksquare$$

%$\qed$
%\end{proof}

Identifying $K_2$ is harder. The ideas of section
\ref{dist_determined_by_moments_sec} give the following
\begin{proposition} \label{maj_k_2}
$K_2(\Maj)\le\omega(1)$
\end{proposition}
\begin{proof}
Consider the distribution of $S_n$ under some $\Q\in\A(n,k,1/2)$.
Obviously, $E_\Q(S_n^l)=E_{\P_{1/2}}(S_n^l)$ for any $l\le k$. The
same holds for $\ol{S_n}$ as it is a linear function of $S_n$.
Therefore, $E_{\Q_n}(\ol{S_n}^l)\ra s_l$ where $s_l=\E(N(0,1)^l)$
is the $l$-th moment of a standard normal distribution. The normal
distribution is determined by its moments. Hence, if
$k(n)\in\omega(1)$ and $\Q_n\in\A(n,k(n),1/2)$ then $\ol{S_n}\ra
N(0,1)$ weakly by proposition \ref{conv_of_moments_prop}. In
particular, $\Q_n(\Maj_n=1)=\Q_n(\ol{S_n}>0)\ra 1/2$.
\end{proof}

In fact for $\Maj$ we can be much more specific.
\begin{theorem} \label{majority_bounds}
There exists a $C>0$ such that for any even $2\le k<n$
\begin{equation}
\frac{C}{\sqrt{k\log k}} \le \max_{Q\in\A(n,k,\frac{1}{2})}|\Q(\Maj_n=1) - \frac{1}{2}| \le
\frac{2\sqrt{2}}{\sqrt{k}}
\end{equation}
And when $\Q_0\in\A(n,n-1,\frac{1}{2})$ is the XOR0 distribution
we have $|\Q_0(\Maj_n=1) - \frac{1}{2}|\ge \frac{1}{3\sqrt{n}}$.
\end{theorem}
The theorem implies that $K_2=\omega(1)$, but is much stronger in that it bounds $\eps^{\Maj}(k,\frac{1}{2})$.

\begin{proof}
The claim about XOR0 is easy to verify directly. The lower bound
comes from a direct construction sketched in the appendix. The upper
bound is actually known in the context of error-correcting codes
\cite[Ch. 9, thm. 23]{MS77} and it appears the proof there also
works in our case. But we point out that a very simple proof of it
can be obtained just by applying theorem
\ref{maximal_mass_difference_thm} and \eqref{rho_middle_value_bound}
to the distribution of $S_n$ and this proof even improves a little
on the constant.
\end{proof}

\subsection{Tribes}

Let $m$ be an integer and let $n=m 2^m$ and let $m(n)$ be its
inverse function. $\Tribes_n$ is the following function: Let the
input bits be divided into $2^m$ sets of size $m$ each, called
\emph{tribes}. Let $y_i$ be the AND of the bits in the $i$-th tribe.
Then $\Tribes_n$ is the OR of the $y_i$'s. Let $S_n=\sum_{0\le i
<2^m} y_i$. Then $\Tribes_n=0$ iff $S_n=0$. Under
$\P_{\frac{1}{2}}$, $S_n\ra Poisson(1)$. It is easily checked that
$\Tribes$ is a sequence of monotone functions with sharp threshold
at $p_c=1/2$ and $\P_{\frac{1}{2}}(\Tribes_n=1)\ra 1-1/e$.

\begin{theorem}
For some $C>0$, \ $C m(n) \le K_1(\Tribes) \le 2 m(n)$
\end{theorem}
\begin{proof}
The proof is similar to that of proposition \ref{maj_k_1}. First,
notice that for $\Q\in\A(n,m(n),p)$ we have $\Q(y_i=1)=p^m$. For
$p<1/2$, a union bound now yields $\Q(\max_{0\le i <2^m} y_i = 1)
\le (2p)^m \ra 0$. If $\Q\in\A(n,2m(n),p)$ then the $y_i$s are
pairwise independent. Using Chebyshev's inequality on $S_n$ yields
the desired result for $p>1/2$.

For the lower bound we use equation
\ref{lower_bound_minimal_n_one_bit_0} to produce a $\Q\in\A(n,C
m(n),p)$ such that the probability of all bits of any tribe are 1 is
0.
\end{proof}

\begin{theorem}
$K_2(\Tribes) \le \omega(m(n))=\omega(\log(n))$
\end{theorem}
\begin{proof}
This is like the proof of \ref{maj_k_2}. There is no need to
normalize $S_n$ as it tends to Poisson(1) as is. Again, we only need
to check that Poisson distribution satisfies Carleman's condition.
\end{proof}

A more refined results, like those for $\Maj$ can be reached using
Theorem \ref{maximal_mass_difference_thm}.

\begin{theorem}
$\eps^{\Tribes_n}(k m(n), 1/2) \le \frac{2}{(k/2)!}$
\end{theorem}

\subsection{$AC^0$ functions} \label{AC0_functions}

$AC^0$ is the class of functions computable by boolean circuits
using Not gates, a polynomial number of AND and OR gates (with
unlimited fan-in) and of bounded depth. $\Tribes$ is a notable
example of an $AC^0$ function of depth 2. Linial and Nisan
(\cite{LN}) conjectured that any boolean circuit of depth $d$ and
size $s$ has $K_2\supset \omega(\log^{d-1} s)$.

We prove a very special case of this conjecture. Let $n=2^{2m}$ and
let the input bits be divided into disjoint sets, $A_i$, consisting
of $m$ bits each. A function is \emph{paired} if it is the OR of AND
gates, each operating on the bits in exactly 2 of the $A_i$'s. A
paired function is, in particular, an $AC^0$ function of depth $2$.

\begin{theorem}
If $f$ is paired then $K_2(f)\le \omega(\log n)$
\end{theorem}

\begin{proof}[proof sketch]
Let $S(f)$ be the number of satisfied AND gates in $f$. The crux of
the proof is to trim $f$ by removing some of the AND gates to
produce a function $f'$, which is (a) very close to $f$ under any
$\omega(\log n)$-wise independent distribution, and (b) $S(f')$
under $\P_p$ tends to a RV which is determined by its moments.
\end{proof}

%The trimming is required mainly for $S$ to converge to a limit
%distribution. Even when it does, this limit does not necessarily
%satisfy Carleman's condition. As an example, consider the AND of 2
%$\Tribes_{m 2^m}$ functions. This may be represented as a paired
%function, by taking all the pairs of tribes of the 2 functions.
%However, $S$ here tends to the product of 2 Poisson(1) distributions
%which does not satisfy Carleman's condition.

\subsection{Majority of majorities}
\label{majority_of_majorities_section}

Let $m$ be an odd integer and let $n=m^2$. $\Maj^2$ is the following
function: divide the input bits into $m$ disjoint sets of size $m$.
Let $y_i$ be the majority of the $i$-th set, then $\Maj^2$ is the
majority of the $y_i$'s.

\begin{theorem}
$K_1(\Maj^2)=2$
\end{theorem}
\begin{proof}
The proof of \ref{maj_k_1} yields $\Q(y_i=1)\le 1/(m (1-2p)^2)$ for
any $\Q\in\A(n,k,p)$, when $p<1/2$. Therefore $E_\Q(\sum_i y_i) \le
1/(1-2p)^2$. The $y_i$'s are not pairwise independent, but Markov's
inequality is enough: $\Q(\Maj^2_n=1)\le 2/(m (1-2p)^2) \ra 0$.
\end{proof}

Notice that this proof applies also to $i$-levels majority,
$\Maj^i$, defined similarly on $m^i$ bits.

\begin{theorem}
$\sqrt{n} \le K_2 \le \omega(\sqrt{n})$
\end{theorem}
\begin{proof}
To show that any function $k\in\omega(\sqrt{n})$ belongs to
$K_2(\Maj^2)$ notice that if $\Q\in\A(n,k,1/2)$ then the
distribution generated on the $y_i$'s belongs to $\A(n,k/m,1/2)$.
Since the $y_i$'s enter majority to produce the output, it is
enough, by theorem \ref{maj_k_2} to have $k/m = \omega(1)$ in order
for the output to tend to 1/2.

To show that $k=m-1$ is not in $K_2$, let $\Q$ be the following
distribution: $\Q$ is $XOR0$ on each $A_i$ and completely
independent on different $A_i$'s. Obviously, $\Q\in\A(n,m-1,1/2)$.
By theorem \ref{majority_bounds}, $\Q(y_i=1)\ge 1/2+ 1/3 \sqrt{m}$
(assume WLOG that $(n+1)/2$ is even). Let $S_n=\sum_{i=0}^{m-1} y_i$
and $\ol{S_n}=(2 S_n-m )/\sqrt{m}$. Since the $y_i$'s are
independent we have that $\ol{S_n}\ra N(a,1)$ where $a=\lim (2
\Q(y_i=1) -1) \sqrt{m} \ge 2/3$. Obviously,
$\Q(\Maj^2=1)=\Q(\ol{S_n}>0)$ is bounded away from $1/2$.
\end{proof}

The surprising fact here is the lower bound of $\sqrt{n}$. First it
shows an example where $K_2$ is much larger then $\omega(K_1)$.
Second, it demonstrates that $L_2$ approximation does not imply
"Sandwich $L_1$" approximation (see section \ref{fourier_sec}).

\subsection{Composition of functions}

$\Maj^2$ is a simple example of composition of functions. What can
we say about compositions in general?

Let $n=ml$ and let $f=g(h_1,..,h_m)$ where the $h_i$'s receive
disjoint sets $A_i$ of $l$ bits each. Assume that $h_i$'s are
balanced with respect to $p_c$ and that $p_c(g)=1/2$.

\begin{theorem} \label{general_composition_theorem}
For $\eps\le\frac{1}{2m}$, $k^f(4m\eps,p_c)\le \sum_i
k^{h_i}(\eps,p_c)$
\end{theorem}
\begin{proof}
$g(y_1,..,y_m)$ can be expressed as a sum of monomials of the form
$\prod y_i \prod (1-y_j)$, each involving all of the $y$'s. We take
the upper and lower "sandwich $L_1$" approximating polynomials of
each $h_i$ (which have degree $k^{h_i}(\eps,p_c)$) and plug the
upper in place of any $y_i$ and one minus the lower in place of any
$(1-y_j)$. This produces a polynomial of degree $k=\sum_i
k^{h_i}(\eps,p_c)$ which bounds $f$ from above. The error of each
monomial, when the distribution is $k$-wise independent is at most
$(1/2+\eps)^m - 1/2^m \le m\eps/2^{m-2}$ for $\eps\le\frac{1}{2m}$.
Summing over the monomials we have an error of no more then $4 m
\eps$. The lower bound is similar.
\end{proof}

This is a very general bound - we did not put any restriction on
$g$, it can even be nonmonotone. For example, $K_2$ for the XOR of
two (or boundedly many) majorities is still $\omega(1)$.

For $0<a<1$, define $\Maj^2_a$ to be the majority of $n^a$
majorities of $n^{1-a}$ bits each. It is easy to see that $K_2\le
\omega(n^{1-a})$. Theorem \ref{general_composition_theorem} gives a
bound of $K_2\le \omega(n^{3a})$. However, using finer properties of
the "sandwich $L_1$" approximating polynomials of $\Maj$, we can do
better.

\begin{theorem}
$K_2(\Maj^2_a)=\omega(n^{\min(a,1-a)})$
\end{theorem}
\begin{proof}
We use the approximating polynomials of the upper $\Maj$ function
(the "$g$") instead of the generic polynomial of theorem
\ref{general_composition_theorem}. These are not only of bounded
degree, but also have small coefficients. This implies that the
resulting polynomial is of degree $O(m)$ and produces an error of
$O(n^{a/2} \eps)$, where $m$ is the degree of the approximating
polynomial of the lower $\Maj$ functions and $\eps$ is their error.
Taking $m=n^a$ gives $\eps=1/\sqrt{n^a}=n^{-a/2}$, as required.
\end{proof}

\subsection{Percolation} \label{percolation_section}

Another very interesting example to consider is that of percolation.
Briefly, \emph{percolation} on a graph $G=(V,E)$ is a distribution
on $\{0,1\}^V$, where we identify the bits with the states
$\{\texttt{open},\texttt{close}\}$. We refer the reader to
\cite{grimmett} for details of the theory of percolation. We denote
the set of all $k$-wise independent percolation with marginal
probability $p$ for every vertex to be open by $\A(G,k,p)$. When $G$
is infinite, we are interested in the probability of existence of an
infinite cluster of open vertices. This event is a boolean function
on infinitely many bits.

\begin{theorem}
For $G=\Z^d$ or $G=\T^d$ (the $d$-ary tree), for any $0<p<1$ and any
$k$ there exist a $\Q\in\A(G,k,p)$ such that there is an infinite
open cluster $\Q$-almost surely, and another such $\Q$ with no
infinite open cluster $\Q$-almost surely.
\end{theorem}

The positive part of this theorem follows from the following 2
theorems about finite versions of percolation. Let $f$ be the
function indicating an open crossing of the $n\times n$ grid.

\begin{theorem}
$2^{\sqrt{\log\log n}}=(\log n)^{1/\sqrt{\log\log n}} \le K_1(f)\le
\omega(\log n)$
\end{theorem}

For the tree case, we need to diverge slightly from the boolean
valued setting. Let $f$ be the number of open paths from the root to
the leaves of $\T^d_n$, the $n$-levels $d$-ary tree.

\begin{theorem}
For any $p$, for $k=C\log n$, there is a $\Q\in\A(\T^d_n,k,p)$ such
that $E_\Q(f)\ge 2$.
\end{theorem}

To this end, we present a way of combining $k$-wise independent
distributions to "amplify" the amount of independence, inspired by
the $(u\ |\ u+v)$ lemma of error-correcting codes. Let $\Z_r$ be the
cyclic group of size $r$. Let $\A^r(n,k)$ be the set of all $k$-wise
independent distributions on vectors $(X_1,\ldots, X_n)\in \Z_r^n$
with each $X_i$ uniform in $\Z_r$. Define $A^r(G,k)$ similarly.

\begin{lemma} \label{X+Y_lemma}
Fix $m\ge 1$. Let $X:=(X_1, \ldots, X_n)\in \A^r(n,k)$. Let
$X^i:=(X^i_j)_{j=1}^n$ be $m$ IID copies of $X$. Let also $Y:=(Y_1,
\ldots, Y_n)\in \A^r(n,2k+1)$ be a vector independent of all the
$X$'s. Then the vector with the following coordinates
\begin{equation}
\begin{matrix}
 X^1_1 + Y_1, &X^1_2 + Y_2, &\ldots, &X^1_n + Y_n,\\
 X^2_1 + Y_1, &X^2_2 + Y_2, &\ldots, &X^2_n + Y_n,\\
 \vdots,&\vdots, &\vdots, &\vdots,\\
 X^m_1 + Y_1, &X^m_2 + Y_2, &\ldots, &X^m_n + Y_n
\end{matrix}
\end{equation}
is in $\A^r(mn,2k+1)$
\end{lemma}

Consult the appendix for a proof of a more general result.

\begin{proof} (Sketch, of theorem)
We build distributions in $A^r(\T^d_n,k)$ such that when we identify
0 with $\texttt{open}$ and the rest with $\texttt{close}$, we get
the desired percolation for $p=1/r$.

The proof goes by induction. For $k=1$ (i.e. no independence) a
suitable distribution is just taking $X_i$ to be identical and $n$
to be large enough.

Assume we have $X\in\A^r(\T^d_n,k)$ such that $E_X(f)\ge 2$. We will
construct a suitable $Z\in\A^r(\T^d_m,2 k +1)$ for $m=n+n^2 k \log
d$. Let $X^i$ be independent copies of $X$ and let
$Y\in\A^q(n,2k+1)$be such that probability of $Y=0$ is maximal,
which is roughly $d^{-nk}$, because there are about $d^n$ RVs in
$Y$. Using lemma \ref{X+Y_lemma} we now assign the RVs in $X^i+Y^i$
to the vertices of $\T^d_m$ such that each is assigned to a subtree
of depth $n$ with root at a level divisible by $n$. Thus, with
probability $d^{-nk}$ we have $Y=0$ and then the open paths form a
Galton-Watson tree with an expectation of $2^{m/n}=2^{1+nk\log d}$.
Thus, the total expectation is $d^{-nk}2^{1+nk\log d}=2$.
\end{proof}

Notice that having an open path from the root to the leaves of
$\T^d_n$ is an $AC^0$ function of depth 2. This is an example of a
rather complicated depth 2 function, very different then the paired
functions considered in section \ref{AC0_functions}. Also, this
function does not exhibit a sharp threshold, thus the need for
different terminology.

\subsection{Fourier transform and $K_2$} \label{fourier_sec}

In this section we consider only the case $p=\frac{1}{2}$. The
quantity $E_\Q(f)$ may be represented using the fourier transform as
$\sum \hat{f}(S) \hat{\Q}(S)$. When we consider $f$ as having values
of $\pm 1$ we have $\sum \hat{f}^2(S) = \sum f^2(S)= 2^n$. Therefore
$\hat{f}^2(S)/2^n$ is a probability measure on all subsets of the
bits, called the \emph{Fourier mass}. Now, a distribution is
$k$-wise independent if and only if all of its Fourier coefficients
of levels between $1$ and $k$ (inclusive) are 0. Therefore, if the
fourier mass of $f$ is supported by the first $k$ levels, then
$\E_\Q(f) = \hat{f}(\emptyset) = \P_{\frac{1}{2}}(f)$. One might
conjecture that if most of the Fourier mass is on the first $k$
levels then $E_\Q(f)$ would be small for all $\Q\in\A(n,k,1/2)$.

In \cite{LMN89} and \cite{H01}, Linial, Mansour and Nisan with an
improvement by H{\aa}stad prove that any $AC^0$ function has its
fourier mass concentrated on the first $O(\log^{d-1} s)$ levels
(where $s$ is the size and $d$ the depth). Had the above conjecture
been true, we would have proved that $K_2=\omega(\log^{d-1} s)$
immediately for any such AC0 function (see section
\ref{AC0_functions}).

However, $\Maj^2$ provides a counterexample for this conjecture, as
its $K_2>\sqrt{n}$ while its fourier mass is concentrated on the
bounded levels, i.e, for any $\eps>0$ there exists $C>0$ such that
all but $\eps$ of the mass is below level $C$. This is because
$\Maj^2$ is a composition of \emph{noise stable} functions and is
therefore noise stable itself (see \cite{BKS}). Of course, $\Maj^2$
is not an $AC^0$ function so this conjecture might still be true in
that domain.

\subsection{Maximal probability that all bits are 1} \label{max_prob_all_bits_1_sec}
In this section we investigate the maximal probability that all the
bits are $1$, i.e, the AND function. At the end of the section two
applications of our bounds are given.

Define $M(n,k,p):=\max_{\Q\in\A(n,k,p)} \Q(\text{All bits are 1})$
then
\begin{theorem} \label{general_upper_bound_all_bits_1_thm}
For even $k$
\begin{equation} \label{general_upper_bound_all_bits_1}
M(n,k,p)\le \frac{p^n}{\P(\bin(n,1-p)\le \frac{k}{2})}
\end{equation}
\end{theorem}
\begin{proof}
Fix $\Q\in\A(n,k,p)$, let $S$ count the number of bits which are 1.
Since $S$ has the same first $k$ moments as a $\bin(n,p)$ the result
follows immediately from theorem \ref{maximal_mass_difference_thm}
and \eqref{rho_last_value} applied to $S$.
\end{proof}

Bound \eqref{general_upper_bound_all_bits_1} is a powerful bound in that it seems to give good results for most ranges of the parameters. Here are some corollaries
\begin{corollary}
\begin{align}
M(n,k,p)&\le
2\sqrt{k}\left(\frac{kp}{2e(1-p)(n-\frac{k}{2})}\right)^{\frac{k}{2}} & \text{For any $n$, even $k$ and $p$} \label{final_upper_bound_all_bits_1}\\
M(n,k,p)&\le 10p^n & \text{$k$ even, $n(1-p)\le\frac{k}{2}$}
\end{align}
\end{corollary}
We add that it is possible to get a result similar to
\eqref{final_upper_bound_all_bits_1} by letting $S$ count the
number of bits which are 1, considering $(S-pn)^k$ and applying
Chebyshev's inequality. Still our approach with theorem \ref{maximal_mass_difference_thm} has the
following advantages. First, it is quite simple as the above proof
of theorem \ref{general_upper_bound_all_bits_1_thm} shows. Second,
it gives \eqref{final_upper_bound_all_bits_1} in all ranges of the parameters $n,k$ and $p$, estimating
$\E(S-pn)^k$ appears to become difficult when $k$ also grows with
$n$, or when $np$ is small. Third, it seems to give slightly
better results, the approach with Chebyshev's inequality apparently does not give
the factor $2$ inside the brackets of
\eqref{final_upper_bound_all_bits_1}.

We are also able to obtain \emph{exact} results, for $k=2,3$. This is done by adapting the closed-form expressions appearing in Boros and Prekopa \cite{BP89} to our settings.
\begin{proposition}
Let $M:=\lfloor(n-1)(1-p)\rfloor$ and $\delta:=\{(n-1)(1-p)\}$ (integer and fractional parts respectively). And also
$N:=\lfloor(n-2)(1-p)\rfloor$ and $\eps:=\{(n-2)(1-p)\}$ then
\begin{align}
M(n,2,p) &= \frac{p}{M+2} +
\frac{\delta^2-\delta(1+p)+p}{(M+1)(M+2)} \\
M(n,3,p) &= \frac{p^2}{N+2} +
\frac{p(\eps^2-\eps(1+p)+p)}{(N+1)(N+2)} = M(n-1,2,p)p \label{exact_result_k_3}
\end{align}
\end{proposition}

For lower bounds on $M(n,k,p)$ and an exact result for small $p$, see the appendix.

We present two applications of our bounds. First a definition, for $q$ a prime-power and $k\ge 2$, a matrix $B\in M_{R\times
n}(GF(q))$ is an OA($n$,$k$,$q$), or an \emph{Orthogonal array of strength $k$
with $q$ levels} (see \cite{MS77} and \cite{HSS99}) if a uniformly
chosen row $(X_1,\ldots, X_n)$ of it has $k$-wise independent entries, each uniform in $GF(q)$. If the rows of $B$ form a linear subspace, then $B$ is called a
\emph{linear orthogonal array} and is referred to by its generator matrix $A\in M_{m,n}$ whose rows are a basis for the rows of $B$. We call $A$ a GOA($n$,$k$,$q$) for short.
%Let $q$ be a prime power and $k\ge 2$, A matrix $B\in M_{M\times
%n}(GF(q))$ is called an \emph{Orthogonal array of strength $k$
%with $q$ levels} (see \cite{MS77} and \cite{HSS99}) if a uniformly
%chosen row $(X_1,\ldots, X_n)$ of it is in $\A^q(n,k)$.
%
%If the rows of $B$ form a linear subspace, then $B$ is called a
%\emph{linear orthogonal array} and is referred to by its generator
%matrix $A\in M_{m,n}$ whose rows are a basis for the rows of $B$.
\begin{enumerate}
\item The bound \eqref{general_upper_bound_all_bits_1} can be used to give another proof of the Rao bound (see \cite{HSS99}) on the minimal size of orthogonal arrays over $GF(q)$. To see this, suppose $B$ is an OA($n$,$k$,$q$) for $k$ even, with $R$ rows. We may assume $B$ contains the all zeroes vector. Consider the distribution $\Q\in\A(n,k,\frac{1}{q})$ obtained by sampling uniformly a row of $B$ and mapping each coordinate to a bit by $0\mapsto 1$, other elements to $0$.
%Choose uniformly a row of $B$ and map each coordinate to a bit using the mapping $0\mapsto 1$, other elements to $0$, to obtain a distribution in $\A(n,k,\frac{1}{q})$.
We have $\Q(\text{All ones vector}) = \frac{1}{R}$, hence
%The probability to get the all ones vector is $\frac{1}{R}$,
by \eqref{general_upper_bound_all_bits_1} we now get $R\ge q^n\P(\bin(n,1-\frac{1}{q})\le \frac{k}{2})$ which is the Rao bound,
%\begin{equation}\label{orthogonal_array_bound}
%R\ge q^n\P(\bin(n,1-\frac{1}{q})\le \frac{k}{2})
%\end{equation}
or using the less refined \eqref{final_upper_bound_all_bits_1} we get $R\ge \left(\frac{2e(q-1)(n-\frac{k}{2})}{k}\right)^{\frac{k}{2}}/2\sqrt{k}$.

We mention in this context that for $q=2$, this lower bound is equal to the bound $m(n,k):=\sum_{i=0}^{\frac{k}{2}} \nchoosek{n}{i}$ which also appeared in \cite{ABI86} (in a more general setting) but we note that for $q=2$ we obtain a somewhat stronger result, the bound \eqref{general_upper_bound_all_bits_1} is in fact an upper bound for the size of any atom of the distribution (by xoring a constant vector), hence for this case we improve slightly the known results by adding that the maximum atom of the distribution is bounded by $\frac{1}{m(n,k)}$, not just the size of the sample space.
\item Let $A$ be a GOA($n$,3,3) with $m$ rows. Since the columns of $A$ are $3$-wise linearly independent, a theorem of Meshulam \cite{M95} implies that $n=O(\frac{3^m}{m})$. Consider the distribution $\Q$ in $\A(n,3,\frac{1}{3})$ obtained by sampling a uniform linear combination of the rows of $A$ and mapping to bits by, say,  $0\mapsto 1$ and $1,2\mapsto 0$. We have $\Q(\text{All ones vector}) \le \frac{1}{3^m} = O(\frac{1}{n\log n})$. In contrast, by equation \eqref{exact_result_k_3} there exist $\Q'\in\A(n,3,\frac{1}{3})$ with $\Q'(\text{All ones vector})=\Omega(\frac{1}{n})$. We deduce that there is an asymptotic difference between what distributions obtained from linear orthogonal arrays (by the above method) and general distributions can achieve. This is interesting since most explicit constructions of $k$-wise independent distributions seem to be based on sampling from linear orthogonal arrays.
\end{enumerate}

\subsection{Minimal probability that all bits are 1} \label{min_prob_all_bits_1_sec}
Define $m(n,k,p):=\min_{\Q\in\A(n,k,p)} \Q(\text{All bits are 1})$.
$m(n,k,p)$ can well be $0$, in fact
\begin{proposition}
When $k<n$ and $p\le\frac{1}{2}$ we have $m(n,k,p) = 0$.
\end{proposition}
Since for $p=\frac{1}{2}$ we can take the XOR0, or XOR1 distributions according to the parity of $n$. And for lower $p$'s we can take the AND of this distribution with a fully independent distribution.

For $p\ge\frac{1}{2}$, define $n_c(k,p):= \min \{n\ |\ m(n,k,p)=0\}$. Our main result of the section are two sided bounds on $n_c(k,p)$
\begin{theorem} \label{bounds_on_nc_thm}
For any $p\ge\frac{1}{2}$
\begin{equation} \label{lower_bound_minimal_n_one_bit_0}
n_c(k,p)\ge \begin{cases}\frac{k}{2(1-p)} + 1 & \text{$k$ even}\\\frac{k+1}{2(1-p)} & \text{$k$ odd}\end{cases}
\end{equation}
and when $1-p=\frac{1}{q}$ for a prime-power $q$ and $C>0$ is a
large constant
\begin{equation} \label{upper_bound_minimal_n_one_bit_0}
n_c(k,p)\le C\frac{k}{1-p}\log(\frac{1}{1-p})
\end{equation}
\end{theorem}
The upper bound is based on the Gilbert-Varshamov bound of error-correcting codes (see \cite{MS77}) and one extra idea. The lower bound poses the main difficulty and for it we need a different aspect of the TCMP. Fix $k$ and $p\ge\frac{1}{2}$ and let $\Q\in \A(n_c(k,p),k,p)$ satisfy $\Q(\text{All ones vector})=0$. Let, as usual, $S$ count the number of $1$'s. $S$ is supported on $[0,n-1]$ and has the first $k$ moments of a $\bin(n,p)$. Theorems in the TCMP show this is only possible if $n_c$ satisfies \eqref{lower_bound_minimal_n_one_bit_0}. The actual verification is technical and involves calculating determinants. Consult the appendix for more details.
%be a distribution always having at least one bit 0. Define $S$ to be the number of bits which are 1 when sampling from $\Q$
%
%
%The upper bound on $n_c$ is not hard to obtain by a construction based on the Gilbert-Varshamov bound of error-correcting codes (see \cite{MS77}). The main difficulty is in proving the lower bound. For this we use a different aspect of the TCMP. Fix $k$ and $p\ge\frac{1}{2}$ and let $\Q\in \A(n_c(k,p),k,p)$ be a distribution always having at least one bit 0. Define $S$ to be the number of bits which are 1 when sampling from $\Q$. We
%note that $S$ has the first $k$ moments as a $\bin(n,p)$ random variable and is supported on the integers between $0$ and $n-1$. The TCMP contains exact criteria for the existence of a distribution with support in $[0,n-1]$ and given moments. Applying these criteria to our problem gives that for $S$ to exist, $n_c$ must specify the above lower bound. The actual application of the criteria involves the calculation of several determinants and is quite technical in nature (and is not presented here). It is interesting that we do not know any other way besides the methods of the TCMP to obtain this result.

As in the previous section, one can deduce from this a result about orthogonal arrays. Suppose $B$ is an OA($n$,$k$,$q$) with the property that each row contains the symbol $0$. If $k$ is even, then necessarily $n\ge 1+\frac{kq}{2}$ and if $k$ is odd then $n\ge\frac{(k+1)q}{2}$.

\section{Open questions}
Below we list some of our main open questions:
\begin{enumerate}
%\item Find an example of a function $f$ (or prove there isn't one)
%whose Fourier mass is concentrated on coefficients of size much
%larger than the $K_2$ of $f$. In the language of approximation it
%means that the function $f$ has a good "sandwich $L_1$"
%approximation of degree $K_2$, but no good $L_2$ approximation of
%that (or slightly larger) degree.
\item Say anything non-trivial about the extremal points of
$\A(n,k,p)$.
\item Is "sandwich $L_1$" approximation stronger then $L_2$
approximation?
\item What is $K_2(\Maj^3)$? What is $K_2(\Maj^i)$?
\item What is $K_2$ for iterated majority of threes?
\end{enumerate}

\section{Appendix}

\subsection{$K_1$ and $K_2$}

It is not a-priori clear whether one of these classes contains the
other. Assume that $\lim \P_{p_c}(f=1)$ exists and denote it by
$\alpha$. We have the following simple result:

\begin{claim}
For any $k\in K_2(f)$, for $p<p_c$ we have $\ol{\lim} \
\eps^f(k,p)\le \alpha$ and for $p>p_c$ we have $\ol{\lim} \
\eps^f(k,p)\le 1-\alpha$
\end{claim}

\begin{proof}
Obviously, both $\max_{\Q\in\A(n,k,p)} \Q(f=1)$ and
$\min_{\Q\in\A(n,k,p)} \Q(f=1)$ are increasing functions of $p$. The
claim now follows immediately from the fact that $\lim
\max_{\Q\in\A(n,k,p)} \Q(f=1) = \lim \min_{\Q\in\A(n,k,p)} \Q(f=1) =
\alpha$.
\end{proof}

So, while we don't know if for $k\in K_2$, $\eps^f(k,p)\ra0$ we do
know that it cannot be too large.

\subsection{Percolation}

Here is the more general result, of which lemma \ref{X+Y_lemma} is a
corollary (put $l=1$).

\begin{lemma} \label{u_u_plus_v_variation_lemma}(combining distributions)
Fix integers $l,m\ge 1$. Suppose for each $1\le i\le m$ we have
random vectors $X^i:=(X^i_1, \ldots, X^i_n)\in \A^r(n,k)$ and
$Y^i:=(Y^i_1, \ldots, Y^i_n)\in \A^r(n,lk+l+k)$. Suppose that the
$X$ vectors are independent among themselves and independent from
the $Y$ vectors, and that the $Y$ vectors are $l$-wise independent
among themselves. Then the vector with the following coordinates
\begin{equation}
\begin{matrix}
 X^1_1 + Y^1_1, &X^1_2 + Y^1_2, &\ldots, &X^1_n + Y^1_n,\\
 X^2_1 + Y^2_1, &X^2_2 + Y^2_2, &\ldots, &X^2_n + Y^2_n,\\
 \vdots,&\vdots, &\vdots, &\vdots,\\
 X^m_1 + Y^m_1, &X^m_2 + Y^m_2, &\ldots, &X^m_n + Y^m_n
\end{matrix}
\end{equation}
is in $\A^r(mn,lk+l+k)$
\end{lemma}
%Before we prove the lemma let us make a few remarks:
%\begin{enumerate}
%\item One special case which may highlight the usefulness of the lemma. Note that taking $l=1$ we can have
%the vectors of $Y$ not independent among themselves, one possibility
%is to take them all to be just the same vector. In this case we have
%just one vector (namely $Y^1$) of $n$ variables which are $2k+1$
%independent, but in the conclusion of the lemma we get a vector of
%$mn$ variables which are $2k+1$ independent (and $m$ may be
%arbitrarily large). In our use, adding the $Y$'s allows us to
%"amplify" the amount of independence while still preserving some of
%the useful properties of the distribution of the $X$'s.
%\item The case when $l=1$, $m=2$, $Y^1=Y^2$ and $X^2=0$ (the zero vector) corresponds to the $(u\ |\ u+v)$ construction of error correcting codes (see \cite[Ch. 2.9, theorem 33]{MS77} or \cite[Sec. 10.3]{HSS99}). This case is not strictly covered by our formulation (since we do not allow the $X^j$ to be zero), but the conclusion still holds and gives $2k+1$ independence.
%\item The lemma is in fact true under more general assumptions, the coordinates of the $X$'s and $Y$'s can be taken to be distributed according to Haar measure in any compact group and, in fact, in any space on which there is a probability distribution $Q$ and an operation "$+$" such that $Q$ is preserved by $x+y$ for any fixed $x$ and for any fixed $y$ (and then the coordinates of $X$ and $Y$ will be distributed according to $Q$).
%\end{enumerate}
\begin{proof}
%Divide the vector $Z$ into $Z^1,\ldots, Z^m$

Call the resulting distribution $Z$, where $Z^i:=(Z^i_1,
Z^i_2,\ldots,Z^i_n)$ and $Z^i_j:=X^i_j + Y^i_j$. Take a set $S$ of
at most $lk+l+k$ variables from the vector $Z$, we need to show they
are independent and uniformly distributed. Suppose that $a^i$ of
them are from $Z^i$ for each $1\le i\le m$, WLOG we can assume that
for each $i$ these are $Z^i_1,\ldots, Z^i_{a^i}$. Consider only the
$i$'s for which $a^i\ge k+1$, since $|S|\le l(k+1)+k$ we can have at
most $l$ such $i$'s, WLOG suppose these are $a^1,\ldots, a^t$ for
$t\le l$. Now, fix some values $c^i_j\in \Z_r$ for $1\le i\le m,\
1\le j\le a_i$, define events $A:=\{Z^i_j = c^i_j\text{ for all
$1\le i\le t$ and $1\le j\le a_i$}\}$ and $B:=\{Z^i_j = c^i_j\text{
for all $t+1\le i\le m$ and $1\le j\le a_i$}\}$. We need to show
that $\P(A,B) = r^{-|S|}$. We start with
\begin{equation}
\begin{split}
\P(A) &= \E\left(\P(A\ |\ (X^i)_{i=1}^t)\right) =\\
&= \E(\P(Y^i_j = c^i_j - X^i_j\text{ for $1\le i\le t$, $1\le j\le a^i$}\ |\ (X^i)_{i=1}^t)) =\\
&= \E(r^{-\sum_{i=1}^t a^i}) = r^{-\sum_{i=1}^t a^i}
\end{split}
\end{equation}
Where the next to last equality follows since the $Y^i_j$ are
uniform and independent from the $X$'s, since they are $l$-wise
independent as vectors (and $l\ge t$), since they are
$(lk+l+k)$-wise independent inside each vector and since $a^i_j\le
|S|\le lk+l+k$.

To finish the lemma we need to show that $\P(B\ |\ A) =
r^{-\sum_{i=t+1}^m a^i}$. We will show something stronger, that in
fact $\P(B\ |\ (X^i)_{i=1}^t,\ (Y_i)_{i=1}^m) = r^{-\sum_{i=t+1}^m
a^i}$. Indeed
\begin{equation}
\begin{split}
\P&(B\ |\ (X^i)_{i=1}^t,\ (Y_i)_{i=1}^m) =\\
&= \P(X^i_j=c^i_j - Y^i_j\text{ for $t+1\le i\le m$, $1\le j\le a^i$}\ |\ (X^i)_{i=1}^t,\ (Y_i)_{i=1}^m) =\\
&= r^{-\sum_{i=t+1}^m a^i}
\end{split}
\end{equation}
Where the last equality follows since the $X^i_j$ for $t+1\le i\le
m$ are uniform, independent from the $Y$'s and from $(X^i)_{i=1}^t$,
since $a^i\le k$ for each $t+1\le i\le m$ by the definition of $t$
and since $(X^i_j)_{j=1}^n$ are $k$-wise independent. This finishes
the proof of the lemma.
\end{proof}

\subsection{The classical moment problem}

Here is the general setup of the classical moment problem
(\cite{Ak65}, \cite{KN77}) leading to the bounds of theorem
$\ref{maximal_mass_difference_thm}$. It is followed by a definition
of the Krawtchouk polynomials and a proof of \eqref{rho_last_value}
and \eqref{rho_middle_value_bound}.

Consider a real sequence $\S := \{s_m\}_{m=0}^k$, with $s_0=1$ (this
last condition is convenient for us in order to use probabilistic
notation, but it is not necessary for the results of the classical
moment problem). Define
\begin{equation}
\A_{\S} = \{\Q\ |\ \text{$\Q$ a probability distribution on $\R$, $s_m
= \E_{\Q} (X^m)$ for $0\le m\le k$}\}
\end{equation}
to be all probability distributions with these first $k$ moments
($X$ is a random variable distributed according to $\Q$).
\begin{definition}
Given $\S = \{s_m\}_{m=0}^k$ with $s_0=1$ and $k$ even, define the
\emph{orthogonal polynomials with respect to $\S$},
$\{P_m\}_{m=0}^{k/2}$ as the unique polynomials with the following
properties:
\begin{enumerate}
\item $P_m$ is a polynomial of degree $m$ with positive leading coefficient.
\item Defining formally a linear operator $T$ from polynomials of degree $k$ to reals by $T(x^i):=s_i$ for $0\le i\le k$ then $T(P_l(x)P_m(x)) = \delta_{l,m}$.
\end{enumerate}
Note that the second condition is the same as requiring $E_{\Q} (P_l(X)P_m(X)) = \delta_{l,m}$ for any $\Q\in\A_{\S}$.
\end{definition}
We remark that these polynomials cannot be defined for degree
larger than $n$ if the sequence $\S$ corresponds to the moments of
an atomic distribution with only $n$ atoms.

Define also the function $\rho_n(x):=(\sum_{m=0}^n P_m^2(x))^{-1}$, then we have the following
\begin{theorem} \cite[2.5.2 and 2.5.4]{Ak65}
For any $x$ and any $\Q_1,\Q_2\in\A_{\S}$
\begin{equation}
|\Q_1(X\le x) - \Q_2(X<x)|\le \rho_{\frac{k}{2}}(x)
\end{equation}
and in particular when $\Q_1=\Q_2$
\begin{equation}
\max_{\Q\in\A_{\S}} \Q(X=x)\le \rho_{\frac{k}{2}}(x)
\end{equation}
\end{theorem}
We remark that in many cases, the theory also has constructions
which achieve these bounds, but we could not use these since in the
cases we needed we required the support of the distribution to be on
integer points. It is possible, however, that a modification of
these constructions can yield a distribution on integer points, this
would be very useful to show the sharpness of the bounds in the
cases we use.

\subsubsection{Krawtchouk polynomials}
In our work we utilize the orthogonal polynomials corresponding to
the moments of the binomial distribution (i.e., when $s_m=\E(X^m)$
where $X\sim\bin(n,p)$). These are the well-known Krawtchouk
polynomials (see \cite{Sz75}). For given $n$ and $p$, the $m$'th
polynomial ($0\le m\le n$) is given by
\begin{equation}
P_m(x)=\nchoosek{n}{m}^{-\frac{1}{2}}\left(p(1-p)\right)^{-\frac{m}{2}}\sum_{j=0}^m
(-1)^{m-j} \nchoosek{n-x}{m-j}\nchoosek{x}{j}p^{m-j}(1-p)^{j}
\end{equation}
where for real $x$ and integer $b\ge1$,
$\nchoosek{x}{b}:=\frac{x(x-1)\cdots(x-b+1)}{b!}$ and
$\nchoosek{x}{0}:=1$.

We note that
\begin{equation}
P_m(n)=\nchoosek{n}{m}^{\frac{1}{2}}\left(\frac{1-p}{p}\right)^{\frac{m}{2}}
\end{equation}
Hence
\begin{equation} \label{rho_m_at_n}
\rho_m(n)=\left(\sum_{j=0}^m
\nchoosek{n}{j}\left(\frac{1-p}{p}\right)^j\right)^{-1} =
\frac{p^n}{\P(\bin(n,1-p)\le m)}
\end{equation}

Furthermore, for $p=\frac{1}{2}$
\begin{equation}
P_m(\frac{n}{2})=\nchoosek{n}{m}^{-\frac{1}{2}}\sum_{j=0}^m
(-1)^{m-j} \nchoosek{n/2}{m-j}\nchoosek{n/2}{j}
\end{equation}
but, as is well known, since the sum is the coefficient of $z^m$
in the power series expansion of
$f(z):=(1+z)^{\frac{n}{2}}(1-z)^{\frac{n}{2}}$ and since
$f(z)=(1-z^2)^{\frac{n}{2}}$ we get by the binomial formula that
\begin{equation} \label{middle_Krawtchouk_value}
P_m(\frac{n}{2}) = \begin{cases}0 & m\text{ odd}\\
\nchoosek{n}{m}^{-\frac{1}{2}}(-1)^{\frac{m}{2}}\nchoosek{n/2}{m/2}
& m\text{ even}\end{cases}
\end{equation}
We then obtain
\begin{lemma} \label{rho_m_at_n_over_2}
For $p=\frac{1}{2}$, even $n$ and even $m\le \frac{n}{2}$
\begin{equation}
\rho_m(\frac{n}{2}) \le \frac{2}{\sqrt{m}}
\end{equation}
\end{lemma}
\begin{proof}
Using \eqref{middle_Krawtchouk_value} we have
\begin{equation} \label{middle_rho_equation}
\rho_m(\frac{n}{2}) = \left(\sum_{j=0}^\frac{m}{2}
\nchoosek{n}{2j}^{-1}\nchoosek{n/2}{j}^2\right)^{-1}
\end{equation}
we recall the well-known inequality that for any integer $a\ge 0$,
$a! = \sqrt{2\pi a}\left(\frac{a}{e}\right)^a e^{\lambda_a}$ where
$\frac{1}{12a+1}\le \lambda_a\le \frac{1}{12a}$. Using this we
notice that
\begin{equation}
\nchoosek{a}{b} = \sqrt{\frac{a}{2\pi
b(a-b)}}\frac{a^a}{(a-b)^{a-b}b^b} e^{\lambda_a - \lambda_{a-b} -
\lambda_b}
\end{equation}
Hence after cancelation
\begin{equation}
\nchoosek{n}{2j}^{-1}\nchoosek{n/2}{j}^2 =
\sqrt{\frac{n}{\pi(n-2j)j}}e^{2\lambda_{n/2}-2\lambda_j-2\lambda_{\frac{n}{2}-j}+\lambda_{n-2j}+\lambda_{2j}-\lambda_n}
\end{equation}
so for $\frac{n}{2},j,(\frac{n}{2}-j)\ge 1$ we get
\begin{equation}
\nchoosek{n}{2j}^{-1}\nchoosek{n/2}{j}^2 \ge
\sqrt{\frac{n}{\pi(n-2j)j}}e^{-\frac{5}{12}} \ge
\sqrt{\frac{1}{8j}}
\end{equation}
Plugging back into \eqref{middle_rho_equation} we get
\begin{equation}
\rho_m(\frac{n}{2}) \le \left(1+\sum_{j=1}^\frac{m}{2}
\sqrt{\frac{1}{8j}}\right)^{-1} \le
\left(1+\frac{1}{\sqrt{2}}\left(\sqrt{\frac{m}{2}}-1\right)\right)^{-1}
\le \frac{2}{\sqrt{m}}
\end{equation}
\end{proof}

\subsection{Majority}
We now continue and give a sketch of the proof of the lower bound
for theorem \ref{majority_bounds}.

\begin{proof} (sketch of lower bound in theorem \ref{majority_bounds})
Fix an odd $n$ and a $2\le k<n$, we would like to construct a distribution $\Q\in\A(n,k,\frac{1}{2})$ such that when we define $S$ to be the number of bits which are $1$ when sampling from $\Q$ then
\begin{equation} \label{lower_bound_for_maj_apdx}
|\Q(S\ge\frac{n+1}{2})-\frac{1}{2}|\ge\frac{C}{\sqrt{k\log k}}
\end{equation}
for some $C>0$. We may assume that $k<c\frac{n}{\log n}$ for some small $c>0$, otherwise the bound follows trivially by taking the distribution XOR0 and using the bound that it satisfies (see theorem \ref{majority_bounds}. Let $M:=C\sqrt{\frac{n}{k\log k}}$ be an integer, the idea of the proof is to construct $\Q$ in such a way that with high probability $S\equiv L$ mod $M$ for some fixed integer $0\le L< M$, and furthermore that on this event $S$ behaves like a $\bin(n,\frac{1}{2})$ random variable conditioned to be $L$ mod $M$. Such an $S$ will satisfy \eqref{lower_bound_for_maj_apdx} for the correct choice of $L$.

To do this, we consider a distribution $\tilde{\Q}$ on $(X_1,\ldots, X_{k+1})\in\Z_M^{k+1}$ satisfying that all the $X_i$ are IID uniform in $\Z_M$ except that $X_k$ is chosen so that their sum is always $L$ modulo $M$. This distribution is of course $k$-wise independent. the required distribution $\Q$ is a distribution on $n$ bits $(Y_1,\ldots, Y_n)$, we create it from the distribution $\tilde{\Q}$ by dividing the $Y$'s into $k+1$ disjoint groups of bits, each $X_i$ is responsible for the value of one of these groups of bits in the following way, when observing the value of the $X$ variable, we sample as uniformly as is possible a string of bits for the $Y$ variables in its group such that their sum modulo $M$ equals the value of the $X$ variable.

The parameters have been chosen in such a way so that the probability that we do not succeed even at one of the $Y$ groups to have the correct sum modulo $M$ is very small. Hence the distribution $\Q$ thus constructed satisfies the required properties.
\end{proof}
\subsection{More bounds on the maximal probability that all bits are 1}
In this section we detail more bounds on the maximial probability that all bits are $1$. Recall $M(n,k,p):=\max_{\Q\in\A(n,k,p)} \Q(\text{All bits are 1})$.

In the main text we have shown
\begin{theorem}
For even $k$
\begin{equation}
M(n,k,p)\le \frac{p^n}{\P(\bin(n,1-p)\le \frac{k}{2})}
\end{equation}
In particular for even $k$
\begin{equation}
M(n,k,p)\le
2\sqrt{k}\left(\frac{kp}{2e(1-p)(n-\frac{k}{2})}\right)^{\frac{k}{2}}
\end{equation}
and for even $k$ and $n(1-p)\le \frac{k}{2}$
\begin{equation}
M(n,k,p)\le 10p^n
\end{equation}
\end{theorem}

We now compliment these with lower bounds on $M(n,k,p)$. Both lower bounds come from well-known constructions of linear error-correcting codes. In both we assume $p=\frac{1}{q}$ for either a prime, or a prime-power, $q$. To get the bounds we first construct the linear code over $GF(q)$, then pass to its dual code, well known to be an orthogonal array. Then sample a line of the orthogonal array uniformly and map to bits using $0\mapsto 1$ and the rest of the elements mapping to $0$. We obtain
\begin{theorem} \label{lower_bound_all_bits_1_thm}
Using the Gilbert-Varshamov bound, for $p=\frac{1}{q}$ with $q$ a prime power
\begin{equation} \label{GV_all_bits_1_bound}
M(n,k,p) \ge p\left(\frac{p(k-1)}{en}\right)^{k-1}
\end{equation}
and using BCH codes, when $p=\frac{1}{q}$ with $q$ a prime (not a prime-power!), $k\equiv 1(\text{mod $q$)}$ and $n+1$ is a power of $q$ then
\begin{equation} \label{BCH_all_bits_1_bound}
M(n,k,p) \ge p\left(\frac{1}{n+1}\right)^{(k-1)(1-p)}
\end{equation}
\end{theorem}
We add that there is a gap in the exponent between these lower bounds and our upper bounds, namely the upper bounds have exponent $\frac{k}{2}$ and the lower bounds have, at best, exponent $(k-1)(1-p)$. We do not know to close this gap but remark that it is also present in the theory of error-correcting codes, for a paper discussing this gap for error-correcting codes and the known results there see \cite{DY04}.

We end this section by remarking on one more exact result, for very small $p$'s
\begin{proposition}
When $p\le \frac{1}{n-1}$
\begin{equation} \label{Poisson_all_bits_1}
M(n,k,p) = p^k
\end{equation}
\end{proposition}
This follows quite simply from a direct construction of the distribution. We start by putting probability $p^k$ on the all ones vector, then all the rest of the probabilities of atoms are determined by being $k$-wise independent with marginal $p$, we check that for this range of $p$'s all these other probabilities are indeed positive. This is the same as the fact that the weight distribution of an MDS code is determined, see \cite{MS77}.
\subsection{More on the minimal probability that all bits are 1}
We remark on the proof of theorem \ref{bounds_on_nc_thm}. The construction of the upper bound on $n_c(k,p)$ goes as follows, we start with an orthogonal array with very good parameters over $GF(q)$ (where now $q=\frac{1}{1-p}$) obtained using the Gilbert-Varshamov bound. We then choose a row uniformly and map each of its coordinates to bits. The mapping is chosen so that in each coordinate exactly one element of $GF(q)$ is mapped to $0$, the rest to $1$, but this element is chosen in a greedy fashion to minimize the chance of not having a $0$ anywhere. When $n$ is large enough compared to $k$ this idea succeeds in giving a distribution with probability $0$ for the all ones vector. This gives the upper bound of the theorem.

As detailed in the main text, the lower bound follows from an existence theorem in the theory of the TCMP, we give this theorem here for easy reference.

Let $X\sim \bin(n,p)$ and define
$s_i:=\E(X^i)$. Define the matrices
\begin{equation}
\begin{split}
A(m)&:=(s_{i+j})_{i,j=0}^m\\
B(m)&:=(s_{i+j+1})_{i,j=0}^m\\
C(m)&:=(s_{i+j})_{i,j=1}^m
\end{split}
\end{equation}
then the classical moment problem states (see \cite{Ak65},\cite{KN77} or \cite{CF91} which contains a survey)
\begin{proposition}
A random variable $S$ with moment sequence $\{s_i\}$ supported on
$[a,b]$ exists if and only if
\begin{enumerate}
\item $k$ is odd and $bA(\frac{k-1}{2})\ge B(\frac{k-1}{2})\ge aA(\frac{k-1}{2})$.
\item $k$ is even, $A(\frac{k}{2})\ge 0$ and
$(a+b)B(\frac{k}{2}-1)\ge abA(\frac{k}{2}-1) + C(\frac{k}{2})$.
\end{enumerate}
where as usual, $A\ge B$ means $A-B\ge 0$ means that $A-B$ is
non-negative definite.
\end{proposition}
\end{document}